\documentclass[11pt,psamsfonts,leqno,oneside,letterpaper]{amsart}
\usepackage[letterpaper,left=1.75truein, top=1truein]{geometry}
\usepackage{amssymb,amsmath,amscd,graphics}
\usepackage{color}

\usepackage[colorlinks,linkcolor=blue,citecolor=blue]{hyperref}
\newtheorem{theorem}{Theorem}[section]
\newtheorem{lemma}[theorem]{Lemma}

\newtheorem{conjecture}[theorem]{Conjecture}
\newtheorem{definition}[theorem]{Definition}

\newtheorem{question}[theorem]{Question}

\def\PSL{\mbox{\rm{PSL}}}

\def\Isom{\mbox{\rm{Isom}}}

\title{Orbifolds and Commensurability}

\author{G. S.  Walsh}
\begin{document}

\maketitle
\begin{abstract}  These are notes based on a series of talks that the author gave at the ``Interactions between hyperbolic geometry and quantum groups" conference held at Columbia University in June of 2009. 
\end{abstract} 

\section{Background on hyperbolic manifolds and orbifolds} 

In understanding manifolds and their commensurability classes, we will find it extremely helpful to employ orbifolds. A manifold is an object locally modeled on open sets in $\mathbb{R}^n$, and an orbifold $\mathcal{O}$ is locally modeled on open sets in $\mathbb{R}^n$ modulo finite groups of Euclidean isometries.  That is, each point $x \in \mathcal{O}$ has a neighborhood modeled on $\tilde U/G$, where $G$ is a finite subgroup of $SO(n)$ and $\tilde U$ is an open ball in $\mathbb{R}^n$.  A {\it geometric orbifold} is the quotient of a simply connected Riemannian manifold $X$ by a discrete subgroup $\Gamma$ of $\Isom(X)$, and we say that $\mathcal{O} = X/\Gamma$ is an $X$-orbifold.  In this case the {\it orbifold fundamental group} is the group $\Gamma$. There are non-geometric orbifolds, but we will only be concerned with geometric ones here.   We will describe the structure of orbifolds through some examples, see \cite{orbifoldbook} for a good description of the details.

\begin{description} 
\item[1] The ``football" is an $S^2$-orbifold which is the quotient of the $S^2$ by the group $\mathbb{Z}/3 \mathbb{Z}$ generated by a rotation of $2 \pi/3$ which fixes the north and south poles.  The {\it ramification locus} of an orbifold $\mathcal{O}$ is the set of  points where any neighborhood is modeled on an open set in $\mathbb{R}^n$ modulo a non-trivial group. The ramification locus in this case is two points, which we label 3, since that is the order of the local group.   The {\it underlying space} $| \mathcal{O}|$ of an orbifold $\mathcal{O}$ is the space obtained from $\mathcal{O}$ by forgetting the orbifold structure, which  is $S^2$ in this case. The two ramification points are both modeled on disks in $\mathbb{R}^2$ modulo a group of rotations, and we call these {\it cone points}. The football is commonly denoted as $S^2(3,3)$.  In general, $M^2(r_1,...r_n)$ is a 2-orbifold with underlying space the 2-manifold $M^2$ and $n$ cone points of orders $r_i$. 

\item[2] A common Euclidean orbifold is $S^2(2,2,2,2)$, which is the quotient of $\mathbb{R}^2$ by the group generated by the translations $(x,y) \rightarrow (x+1,y), (x,y) \rightarrow (x, y+1)$ and rotation by $\pi$ about $(0, 1/2)$.  Note that this will generate rotations by $\pi$ in all half-integer lattice points. Exercise: What is a fundamental domain?  Find a rotation of order 4 on $\mathbb{R}^2$ such that the group generated by this rotation and the above generators yields the orbifold $S^2(2,4,4)$.

\end{description}
An {\it orbifold covering} $f: Q' \rightarrow Q$ is a continuous map between the between the underlying spaces $|Q'| \rightarrow |Q|$.  We further require that if a point $x \in |Q|$ has a neighborhood whose orbifold structure is $U = \tilde U/G$ then each component $V_i$ of $f^{-1}(U)$ is isomorphic to $\tilde U/G_i$ where $G_i <G$ and $f|_{V_i}:V_i \rightarrow U$ is $\tilde U/G_i \rightarrow \tilde U/G$.   See \cite[Section 2.3]{orbifoldbook} for elaboration.  Example 1 above gives a cover of $S^2(3,3)$ by $S^2$, and example 2 gives a cover of $S^2(2,4,4)$ by $S^2(2,2,2,2)$.  Note that $S^2(2,4,4)$ also covers itself. 

\begin{description} 
\item[3] 
We regard hyperbolic 2-space ${\bf H}^2$ as the upper-half-space $\lbrace z| Im(z)>0 \rbrace$ of the complex plane and declare $\Isom^+({\bf H}^2)$ to be $\PSL(2,\mathbb{R})$, where the matrix $\bigl( \begin{smallmatrix} 
  a & b\\
  c & d 
\end{smallmatrix} \bigr)$ acts by $z \rightarrow \frac{az+b}{cz+d}$.  Then the group generated by $\bigl( \begin{smallmatrix} 
  2 & 1\\
  1 & 1  
\end{smallmatrix} \bigr)$ and $
\bigl( \begin{smallmatrix} 
  1 & 1\\
  1 & 2 
\end{smallmatrix} \bigr)$ is a free group of rank 2.  The quotient of ${\bf H}^2$ by this group is a punctured torus. If we add a generator that rotates by $\pi$ about $i$, $\bigl( \begin{smallmatrix} 
  0 & 1\\
  -1 & 0 
\end{smallmatrix} \bigr)$, the quotient orbifold is $S^2(2,2,2,\infty)$ (the $\infty$ denotes a cusp).  This is a quotient of the punctured torus by an involution which fixes three points on the punctured torus and takes the puncture to itself.  It is perhaps easier to understand the quotient by considering the induced action on 
the unit disk. Use the transformation $U(z) = \frac{zi+1}{z+i}$ to take the upper-half space to the unit disk. Then the rotation by $\pi$ about $i$ becomes rotation by $\pi$ about the origin.   This illustrates $S^2 (2,2,2,\infty)$ as a hyperbolic 2-orbifold. 
\end{description} 

Here we will mainly consider hyperbolic 3-orbifolds.  These are the quotient of ${\bf H}^3$ (regarded as $\lbrace (z,t)| z \in \mathbb{C}, t \in \mathbb{R}, t >0 \rbrace$) by a discrete finitely generated subgroup of $\Isom^+({\bf H}^3) \cong \PSL(2, \mathbb{C})$.  As above, $\PSL(2, \mathbb{C})$ acts on the complex plane $\lbrace (z,0), z \in \mathbb{C} \rbrace $ such that $\bigl( \begin{smallmatrix} 
  a & b\\
  c & d 
\end{smallmatrix} \bigr)$ acts by $z \rightarrow \frac{az+b}{cz+d}$.  The action on ${\bf H}^3$ is by Poincar\'e extension, where hemispheres perpendicular to the complex plane are mapped to hemispheres perpendicular to the complex plane. (See  \cite[Chapter 1]{outercircles}.) Volume will be discussed in other lectures here, but for our purposes all subgroups of $\PSL(2, \mathbb{C})$ will be discrete and have finite co-volume, meaning that the quotient is a finite-volume hyperbolic orbifold. A discrete subgroup of $\PSL(2, \mathbb{C})$ is called a Kleinian group. Let $\Gamma_1$ and $\Gamma_2$ be two finite co-volume Kleinian groups.  Mostow-Prasad rigidity states that if $ \Gamma_1$ and $ \Gamma_2$ are isomorphic, then they are conjugate. This means that if $\mathcal{O}_1$ and $\mathcal{O}_2$ are finite volume hyperbolic orbifolds ${\bf H}^3/\Gamma_1$ and ${\bf H}^3/\Gamma_2$, then any isomorphism of their fundamental groups is determined by a unique isometry of the hyperbolic orbifolds.  Therefore, invariants which depend only on the conjugacy class of the representation of the fundamental group into $\PSL(2, \mathbb{C})$ are topological invariants.  We will describe two here, the trace field and the cusp field.  Others are discussed in Section \ref{day2}. 

Let ${\bf H}^3/\Gamma$ be a hyperbolic orbifold.  Then the field generated by the traces of elements in $\Gamma$ is invariant under conjugation of $\Gamma$ and this is the {\it trace field} of $\Gamma$.  If there is a subgroup of $\Gamma$  which fixes some point on the sphere at infinity,  $S^2_{\infty} = \mathbb{C} \cup \infty$, then we can conjugate so that this subgroup fixes $\infty$.  Then the quotient of the plane $\lbrace (z, \epsilon)| z \in \mathbb{C} \rbrace$ in ${\bf H}^3$ will be a Euclidean 2-orbifold $\mathcal{Q}$ in ${\bf H}^3/ \Gamma$ for some $\epsilon$, and we call the quotient a {\it cusp cross-section}. The quotient of $\lbrace (z,t)| t >\epsilon \rbrace$ is $\mathcal{Q}^2 \times \mathbb{R}^+$ and we call this a {\it cusp}.  If the cusp cross-section is a torus, we can conjugate so that the generators are $\bigl( \begin{smallmatrix} 
  1 & 1\\
  0 & 1 
\end{smallmatrix} \bigr)$  and $\bigl( \begin{smallmatrix} 
  1 & g\\
  0 & 1 
\end{smallmatrix} \bigr)$.  The {\it cusp field} of the cusp is $\mathbb{Q}(g)$.  Note that $g$ coincides with the shape of the torus in the identification of the Teichmuller space of $T^2$ with upper half-space. If the cusp-cross section is a compact Euclidean orbifold which is not a manifold,  then we take a cover of ${\bf H}^3/\Gamma$ where this cusp cross-section lifts to a torus and take the cusp field of some cusp in the pre-image.  Note that the cusp field of a particular cusp is left unchanged  after taking finite covers.  

\begin{description}

\item[4] The figure-8 knot complement can be realized as a hyperbolic manifold ${\bf H}^3/\Gamma_k$ where  $\Gamma_k$ is generated by $\bigl( \begin{smallmatrix} 
  1 & 1\\
  0 & 1  
\end{smallmatrix} \bigr)$  and $\bigl( \begin{smallmatrix} 
  1 & 0\\
 -w & 1
\end{smallmatrix} \bigr)$, where $w = \frac{-1 + \sqrt{3}i}{2}$.  (This is not obvious, see Thurston's notes \cite[Section 4.3]{thurstonnotes}.)   There is a subgroup $\mathbb{Z} \oplus \mathbb{Z}$ which fixes infinity generated by $\bigl( \begin{smallmatrix} 
  1 & 1\\
  0 & 1
\end{smallmatrix} \bigr)$ and $\bigl( \begin{smallmatrix} 
  1 & g\\
  0 & 1
\end{smallmatrix} \bigr)$ where $g = 2 + 4\sqrt{-3}$. Then the cusp field is $\mathbb{Q}(\sqrt{-3})$ and the trace field is also $\mathbb{Q}(\sqrt{-3})$.  For most, (but not all!) hyperbolic knots in $S^3$ the trace field is the same as the cusp field of the unique cusp. 
\end{description}

\begin{question}   (See \cite{NR}.) What are necessary or sufficient conditions for a hyperbolic knot complement to have its trace field strictly larger than its cusp field?  Currently, there are four knots known that have this property, see \cite{GHH}. 
\end{question}

 A knot $K$ in $S^3$ is {\it strongly invertible} if there is an order two involution $t$ of $(S^3,K)$ such that the fixed point set of $t$ intersects the knot twice.  A strong involution induces an order two involution on the knot complement where the fixed point set intersects each cusp cross-section 4 times.  The figure-8 knot is strongly invertible and we can take the quotient of the figure-8 knot complement by the strong inversion to obtain an hyperbolic 3-orbifold.  This corresponds to adjoining the element $ \tau= \bigl( \begin{smallmatrix} 
  i & 0\\
  0 & -i
\end{smallmatrix} \bigr)$ to the group $\Gamma_k$ in example 4 above.  Note the action of $\tau$ on the cusp subgroup. Also, $\tau$ takes each generator of the knot group to its inverse. 

Let $K$ be a strongly invertible hyperbolic knot. The quotient of $S^3$ by a strong inversion is $S^3$ with an unknotted circle labeled two. The quotient  of $N(K)$, a regular neighborhood of the knot, is a ball with two unknotted arcs labeled 2.  Therefore, the quotient of $S^3 \setminus N(K)$ by the inversion has underlying space a ball and two arcs of the ramification locus which are both labeled two. All of the information about this orbifold is contained in the ramification locus, since its underlying space is topologically trivial.  

\section{Commensurability} \label{day2}

\begin{definition} Two orbifolds are \emph{commensurable} if they admit homeomorphic finite-sheeted covers.
\end{definition} 

We claim that the relation $\sim$ on $n$-dimensional orbifolds defined by commensurability is an equivalence relation.  Indeed, clearly $\mathcal{O} \sim \mathcal{O}$, and if $\mathcal{O} \sim \mathcal{P}$, $\mathcal{P} \sim \mathcal{O}$.   Now assume that $\mathcal{O} \sim \mathcal{S}$ with a common finite index cover $X$, and $\mathcal{S} \sim \mathcal{P}$ with common finite index cover $Y$, so that $X$ and $Y$ both cover $S$.   To show $\mathcal{O} \sim \mathcal{P}$, we will construct a finite index orbifold cover of $X$ and $Y$.  This is the {\it orbifold fibered product} $X \times_S Y$:
\eject
\centerline{$X \times_S Y$}
\centerline{  $\swarrow$ \hskip .3 in $\searrow$}
\centerline{$X$ \hskip .6 in $Y$}
\centerline{$\swarrow$ \hskip .1 in  $\searrow$ \hskip .2 in $\swarrow$  \hskip .1 in $\searrow$ }
\centerline{$\mathcal{O}$\hskip .5 in $\mathcal{S}$\hskip .5 in  $\mathcal{P}$}

Let $f: X \rightarrow S$ and $g:Y \rightarrow S$ be orbifold covering maps.  Then each point $w \in S$ has a neighborhood $U_w = \tilde U / G$ such that a pre-image of $U_w$ in $X$, respectively $Y$, has the form $\tilde U/G_1$, respectively $\tilde U/G_2$  for some $G_i<G$.  It is helpful to think of $\tilde U/G$ as the set of orbits $\lbrace Gy | y \in \tilde U \rbrace$. We define $f_g: \tilde U \rightarrow \tilde U/G_1 \times \tilde U/G_2$ by $f_g(y) = (G_1gy, G_2y)$.  This map factors through $\tilde U/(g^{-1}G_1g \cap G_2)$ since $y$ and $xy$ where $x \in g^{-1}G_1g \cap G_2$ have the same image. We get maps from $\tilde U/(g^{-1}G_1g  \cap G_2)$ to $\tilde U/G_i$ by composing $f_g$ with the projections.    Thus we define the fiber product of $\tilde U/G_1$ and $\tilde U/G_2$ over $\tilde U/G$ as the disjoint union of the orbifolds $\tilde U/(g^{-1}G_1g \cap G_2)$ where $g$ is taken over representatives of the double cosets of $G_1 \backslash G \slash G_2$.  (The  map $f_{g_1gg_2}$ differs from $f_g$ only by the action of $G$ on $\tilde U$.) We can patch these neighborhoods together to get an orbifold which covers $X$ and $Y$.  See \cite[Section 13.2]{thurstonnotes} for an example of why this is the correct way to extend the definition of the fibered product for manifolds.  This exhibits  $X \times_S Y$ as an orbifold cover of  $\mathcal{O}$, $ \mathcal{S}$, and $\mathcal{P}$ so  $\sim$ is transitive which shows the claim.

Aside from the intrinsic reason that this equivalence relation is a method of organizing manifolds, there are many reasons to study commensurability classes.  In particular, there are lots of properties which are preserved by commensurability.  For example, if a manifold $X$ is virtually fibered, and $Y$ is commensurable with $X$, then $Y$ is also virtually fibered. This is because a cover of a manifold which fibers over the circle fibers over the circle.  Similarly, the properties of being virtually Haken, being virtually large, and containing an immersed geodesic surface are all properties of the commensurability class.  In addition, if we restrict ourselves to 3-manifolds which admit a finite-volume geometry, 
this equivalence relation on 3-manifolds preserves the geometric type.  Thus we can regard commensurability classes of geometric 3-manifolds as a refinement of geometrization. This is more useful for some geometries than for others. For example, all spherical orbifolds are commensurable, but there are infinitely many commensurability classes of hyperbolic orbifolds.  

Commensurability classes are particularly relevant for the study of finite-volume hyperbolic 3-manifolds and orbifolds, where classification has often been centered around notions of volume. Complementing this, commensurability classes are transverse to volume. Selberg's lemma states that a finitely generated subgroup of $GL(n,\mathbb{C})$ has a torsion-free subgroup of finite index.  Therefore, all hyperbolic orbifolds are finitely covered by manifolds and there are manifolds in every commensurability class. We define the volume of a hyperbolic orbifold to be $1/d$ the volume of a $d$-fold cover which is a manifold.  By rigidity, volume is a topological invariant and volumes of commensurable manifolds are rationally related.  There has also been recent progress in understanding the commensurator of infinite volume hyperbolic manifolds \cite{LLR}.

Henceforth, we will mainly be concerned with commensurability classes of finite volume hyperbolic 3-orbifolds.    In this case, we can characterize commensurability of hyperbolic orbifolds by commensurability of subgroups of $\PSL(2, \mathbb{C})$. By Mostow-Prasad rigidity, $M_1 \cong {\bf H}^3/\Gamma_1$ and $M_2 \cong {\bf H}^3/\Gamma_2$ are commensurable if and only if  $\Gamma_1$ and a conjugate of $\Gamma_2$ share a finite index subgroup.  In this case we say that the groups $\Gamma_1$ and $\Gamma_2$ are {\it commensurable in the wide sense} in $\PSL(2,\mathbb{C})$.  

Let $\Gamma$ be a finitely generated Kleinian group with finite co-volume.  Then by rigidity, the trace field of $\Gamma$ is a topological invariant of the hyperbolic orbifold ${\bf H}^3/\Gamma$, as mentioned in Section 1.  It is also a finite degree extension of $\mathbb{Q}$ \cite[Theorem 3.1.2]{MR}.  It is not, however, an invariant of the commensurability class.  In the above example of the figure-8 knot complement and the orbifold which is the quotient by the strong inversion, it can be seen that the trace field of the quotient orbifold contains $i$ where the trace field of the figure-8 knot complement is $\mathbb{Q}(\sqrt{-3})$, which does not contain $i$.  However, the {\it invariant trace field } $
k\Gamma= \mathbb{Q}(tr(\gamma^2)| \gamma \in \Gamma)$ is an invariant of the commensurability class \cite[p 117]{MR}.  Often we will restrict to hyperbolic knot complements and in this case the trace field is the same as the invariant trace field \cite{Re}.  

An equivalent way to define the invariant trace field for a cusped finite-volume hyperbolic manifold with an ideal triangulation is the field generated by the cross-ratios of all ideal tetrahedra, or the shapes of all the tetrahedra \cite[Section 5.5]{MR}.  This suggests that the cusp field is a subfield of the invariant trace field and indeed this is the case. 

Probably the most useful invariant of a commensurability class is the commensurator.  Two subgroups $A$ and $B$ of a group are commensurable if their intersection is finite index in each.  The {\it commensurator}  of a Kleinian group $\Gamma$  is  
$$C^+(\Gamma) = \lbrace g \in \
\PSL(2, \mathbb{C}) | \ g \Gamma g^{-1} \text{and}\  \Gamma \ \text{are commensurable} \rbrace .$$ 

The geometry of the commensurator is dramatically different for arithmetic and non-arithmetic Kleinian groups. 

 A Kleinian group $\Gamma$ is {\it arithmetic} if it is commensurable with a $k$-embedding into $M_2(\mathbb{C})$ of the group of norm 1 elements of an order of a quaternion algebra $A$ over $k$ where $A$ ramifies at all real places and where $k$ is a number field with one complex place.  (See \cite[Definition 8.2.1]{MR} and the discussion therein.)  A useful characterization is that a non-cocompact Kleinian group $\Gamma$ is arithmetic if and only if a conjugate of $\Gamma$ is commensurable with a Bianchi group. These are the groups $\PSL(2, O_d)$, where $O_d$ is the ring of integers in $Q(\sqrt{-d})$.  We will say that the orbifold ${\bf H}^3/\Gamma$ is arithmetic when $\Gamma$ is arithmetic.  For example, the complement of the figure-8 knot described above is arithmetic and furthermore this is the only arithmetic knot complement in $S^3$ \cite{Re1}.

Let $\Gamma$ be a finitely generated Kleinian group with finite co-volume.  Margulis's theorem (which holds in more generality) implies that either $\Gamma$ is finite index in $C^+(\Gamma)$, or $C^+(\Gamma)$ is dense in $\PSL(2, \mathbb{C})$. Furthermore, $C^+(\Gamma)$ is dense exactly when $\Gamma$ is arithmetic.

\begin{lemma} 
For finite co-volume Kleinian groups $\Gamma_1$ and $\Gamma_2$, the orbifolds ${\bf H}^3/\Gamma_1$ and ${\bf H}^3/\Gamma_2$ are commensurable if and only if the commensurators of $\Gamma_1$ and $\Gamma_2$ are conjugate in $\PSL(2,\mathbb{C})$.
\end{lemma} 

The fact that commensurable orbifolds have conjugate commensurators follows directly from Mostow-Prasad rigidity and the definition.  For non-arithmetic orbifolds, Margulis's theorem says that $\Gamma_1$ and $\Gamma_2$ are finite index in their commensurators, so if their commensurators are conjugate, the orbifolds ${\bf H}^3/\Gamma_1$ and ${\bf H}^3/\Gamma_2$ are commensurable.  For the arithmetic case, the commensurators of $\Gamma_1$ and $\Gamma_2$ are the set of invertible elements in the quaternion algebra defined by the trace field. (See \cite{MR}.) Therefore we may assume, after conjugation, that $\Gamma_1$ and $\Gamma_2$ are contained in two orders $O_1$ and $O_2$ in the same quaternion algebra. Furthermore, the intersection of two orders is an order so we may assume that $O_1 \subset O_2$.  Since we are considering groups up to commensurability and conjugation, the result follows  from the fact that if $O_1'$ and $O_2'$ are the groups of elements of norm 1 in each order, then $O_1'$ is finite index in $O_2'$. \qed

Thus for a non-arithmetic finite-covolume Kleinian group $\Gamma$,  the commensurator $C^+(\Gamma)$ is the maximal element in the commensurability class of a Kleinian group $\Gamma$.  In the non-arithmetic case, this corresponds to an orbifold ${\bf H}^3/C^+(\Gamma)$ which is the minimal element in the commensurability class of $\mathcal{O} = {\bf H}^3/\Gamma$.  Namely every orbifold commensurable with $\mathcal{O}$ finitely covers ${\bf H}^3/C^+(\Gamma)$.  Thus two non-arithmetic finite-volume hyperbolic orbifolds have a finite sheeted-cover exactly when they finitely cover a common orbifold. 

\section{Hidden Symmetries}

Symmetries of a hyperbolic manifold $M$ and symmetries between finite covers of $M$ will play a very important role in understanding the commensurability class of $M$. Recall that for a Kleinian group $\Gamma$, the {\it normalizer} of $\Gamma$ is 

$$N^+(\Gamma) = \lbrace g \in \PSL(2, \mathbb{C})| g \Gamma g^{-1} = \Gamma \rbrace .$$

Any self-isomorphism of a finite volume hyperbolic manifold $M = {\bf H}^3/\Gamma$ yields an automorphism of $\Gamma$ which by Mostow-Prasad rigidity is realized by conjugation. $\Gamma \leq N^+(\Gamma)$ and these conjugations are realized by a change of base point. Therefore, $Isom^+({\bf H}^3/\Gamma) \cong N^+(\Gamma)/\Gamma$. 

The hyperbolic orbifold ${\bf H}^3/N^+(\Gamma)$ may or may not be minimal in the commensurability class of ${\bf H}^3/\Gamma$. Clearly, $N^+(\Gamma) \leq C^+(\Gamma).$

\begin{definition} Let $\Gamma$ be a finite co-volume Kleinian group.  If $N^+(\Gamma)$ is strictly smaller than $C^+(\Gamma)$, then $\Gamma$ (and ${\bf H}^3/\Gamma$) are said to have \emph{hidden symmetries}. 
\end{definition} 

The elements of the commensurator correspond to isometries between finite-sheeted covers of ${\bf H}^3/\Gamma$.  Indeed, if $g \Gamma g^{-1} \cap \Gamma$ is a finite-index subgroup of $\Gamma$, then conjugation by $g^{-1}$ gives an isomorphism between the subgroups $g \Gamma g^{-1} \cap \Gamma$ and $\Gamma \cap g^{-1} \Gamma g$, which is an isometry between the corresponding finite-sheeted covers.  If $\phi: {\bf H}^3/\Gamma_1 \rightarrow {\bf H}^3/\Gamma_2$ is an isometry between finite-sheeted covers, then it is realized by conjugation by $g$, for some $g \in \PSL(2, \mathbb{C})$  by Mostow-Prasad rigidity.  Therefore, $\Gamma$ and $g \Gamma g^{-1}$ are commensurable. 

Any isometry of ${\bf H}^3/\Gamma$ will permute the finite-sheeted covers of ${\bf H}^3/\Gamma$. Thus, if the commensurator of $\Gamma$ is strictly larger than the normalizer, there is an isometry between finite-sheeted covers of ${\bf H}^3/\Gamma$ which is not realized by an isometry of ${\bf H}^3/\Gamma$.  Hence the term ``hidden symmetries".  Also, each element of the commensurator yields an element of the normalizer of a finite index subgroup $\Gamma'$ of $\Gamma$.  In the non-arithmetic case, by Margulis's theorem above, $\Gamma$ is finite-index in $C^+(\Gamma)$.  Therefore, for every $g \in C^+(\Gamma)$, there is a $k$ such that $g^k \in \Gamma$.  Hence $\Gamma \cap  g \Gamma g^{-1} \cap  g^2 \Gamma g^{-2} ...\cap  g^k \Gamma g^{-k}$ is normalized by $g$ for some $k$.  In the arithmetic case, there are infinitely many hidden symmetries, see \cite{FW}. In the non-arithmetic case, there is a finite-sheeted cover of ${\bf H}^3/\Gamma$ which exhibits all the hidden symmetries of ${\bf H}^3/\Gamma$, namely ${\bf H}^3/\Gamma'$, where $\Gamma'$ is the intersection of all the conjugates of $\Gamma$ in $C^+(\Gamma)$. 

Goodman, Heard and Hodgson \cite{GHH} have recently obtained a characterization of the commensurator of cusped hyperbolic manifolds. Their characterization is as follows.  A {\it horoball packing} is a collection of disjoint horoballs in ${\bf H}^3$.  A {\it cusp neighborhood} is any neighborhood of the cusp that lifts to a horoball packing.  They show \cite[Lemma 2.3]{GHH} that two hyperbolic cusped orbifolds of finite volume cover a common orbifold if and only if they admit choices of cusp neighborhoods lifting to isometric horoball packings. Furthermore, the commensurator of a finite-covolume non-arithmetic Kleinian group $\Gamma$ is the maximal symmetry group of a horoball packing, amongst those which are lifts of cusp neighborhoods in ${\bf H}^3/\Gamma$.  Equivalently, this is the maximal symmetry group of a tiling of ${\bf H}^3$ amongst those which are obtained by lifting canonical cell decompositions of ${\bf H}^3/\Gamma$.  Using this, they are able to compute the commensurators for a large number of cusped hyperbolic orbifolds and to detect the presence of hidden symmetries, by comparing the commensurator with the normalizer. In particular, their computations yield that out of all hyperbolic knots up to 12 crossings, there are only two whose complements admit hidden symmetries.  These are the two dodecahedral knots of  Aitchison and  Rubinstein \cite{AR}.  These are the only two knots whose complements decompose into two regular ideal hyperbolic dodecahedra.  They are commensurable, and one is fibered while the other is not.  

A non-arithmetic hyperbolic 3-orbifold ${\bf H}^3/\Gamma$ admits hidden symmetries if and only if it non-normally covers a finite orbifold, since in this case the orbifold ${\bf H}^3/N^+(\Gamma)$ is not the minimal element in the commensurability class.  We can understand the commensurator orbifold ${\bf H}^3/C^+(\Gamma)$ in part by looking at its cusp cross section, which is a Euclidean 2-orbifold.  The symmetry group of a hyperbolic knot complement is either cyclic or dihedral, since it is finite, and any symmetry of the complement will take a minimal genus Seifert surface (along with the canonical longitude) to a minimal genus Seifert surface.  Therefore, any normal covering of an orbifold by a knot complement has a torus or a $S^2(2,2,2,2)$ cusp cross-section.  Conversely, if a knot complement covers an orbifold with a torus or $S^2(2,2,2,2)$ cusp cross-section, the covering is normal, \cite{Re1, GAW}.  The orientable Euclidean orbifolds which are not a torus or $S^2(2,2,2,2)$ are $S^2(2,4,4)$, $S^2(3,3,3)$ and $S^2(2,3,6)$. These are {\it rigid}, meaning that their moduli spaces are trivial. Thus Neumann and Reid show:

\begin{lemma} \cite[Proposition 9.1]{NR} The following are equivalent for a hyperbolic knot complement ${\bf H}^3/\Gamma$ which is not the figure-8 knot complement. 
\begin{enumerate} 
\item  The knot complement admits hidden symmetries.  \item ${\bf H}^3/C^+(\Gamma)$ has a rigid Euclidean cusp cross-section. \item The knot complement non-normally covers some orbifold. \end{enumerate}  
\end{lemma}
Since a rigid Euclidean 2-orbifold cannot be deformed, the cusp field of an orbifold with one cusp which has a rigid cusp cross-section is either ${\bf Q}(\sqrt{-3})$ (when the cusp cross-section is $S^2(3,3,3)$ or $S^2(2,3,6)$ or ${\bf Q}(i)$ (when the cusp cross-section is $S^2(2,4,4)$.) Furthermore, any hyperbolic orbifold which covers an orbifold with one cusp with a rigid cusp-cross section has the same cusp field. Therefore, a knot complement which admits hidden symmetries has cusp field ${\bf Q}(\sqrt{-3})$ or ${\bf Q}(i)$.  Although it is not known exactly how cusp fields are distributed amongst hyperbolic knot complements, it is suspected that these cusp fields are not any more prevalent than others. 
This, along with the experimental results above, suggests that hyperbolic knots with hidden symmetries are extremely rare. 
 
\section{knot complements} 
  
Although understanding commensurability classes for arbitrary hyperbolic 3-manifolds is quite difficult, the situation for hyperbolic knot complements appears to be much easier.   If the knot complements $S^3 \setminus K$ and $S^3 \setminus K'$ are commensurable, we say that the knots are commensurable and write $K \sim K'$. There are finitely many knots $K'$ with $K \sim K'$  when $K$ (equivalently $K'$) does not admit hidden symmetries.  In this case, the commensurator quotient orbifold will have a flexible cusp, and any filling of the knot complement will cover a filling of the commensurator quotient.  In particular, the $S^3$ filling of the knot complement will cover a filling of the orbifold which necessarily has finite fundamental group.  There are finitely many fillings of a hyperbolic orbifold which yield orbifolds with finite fundamental group, hence there can be only finitely many knot complements in the commensurability class of a knot which does not admit hidden symmetries.  As of this writing, it is not known if there can be infinitely many knots in the commensurability class of a knot which admits hidden symmetries, or even if the commensurability class of the two dodecahedral knots contains infinitely many knot complements, or if there are infinitely many different commensurability classes with hidden symmetries. 

There are several known ways for hyperbolic knots to be commensurable.  The first is the classic (but rare) situation when a knot complement $S^3 \setminus K$ admits a Dehn filling such that the resulting manifold $S^3 \setminus N(K)(r)$ is a lens space.   In this case the universal cover of the filled manifold is $S^3$, with cyclic covering group.  Since a cyclic cover of a knot complement is unique and has one cusp, the pre-image of $K$ in the covering $S^3 \rightarrow S^3 \setminus N(K)(r)$ has one component $K'$ and this induces a covering $S^3 \setminus K' \rightarrow S^3 \setminus K$.  Thus, if the covering and the knot are both non-trivial, then $K \sim K'$. In \cite{GAW} it is proven that a knot complement covers another knot complement if and only if the covered knot admits a cyclic surgery.  By the cyclic surgery theorem \cite{CGLS}, a hyperbolic knot can admit at most three cyclic surgeries.  By volume considerations, different cyclic surgeries will correspond to different knot complements in the commensurability class by the construction above.  Thus, we can only obtain three knot complements in a commensurability class using this method.   There is also a family of commensurable pairs of knot complements given by Walter Neumann. (See \cite{GAW} for a description.) Each pair turns out to cyclically cover a common orbifold.  In addition, there is one example of a pair of knot complements that non-normally cover a common orbifold,   the pair of dodecahedral knots discussed above.  This led to the following conjecture.  Let $K$ be a hyperbolic knot and let $C(K)$ denote the set of knots commensurable with $K$. 

\begin{conjecture} \cite{RW} \label{3conj}
$|C(K)| \leq 3$
\end{conjecture}

This conjecture has been verified in a number of cases, including all hyperbolic two-bridge knot complements \cite{RW}, $(-2,3,n)$-pretzel knot complements \cite{MM}, and for an infinite family of hyperbolic knot complements constructed by Hoffman \cite{Hof} each of which has exactly three knots in its commensurability class. In each of these cases, the knots in question are shown not to admit hidden symmetries by proving that neither ${\bf Q}(i)$ nor ${\bf Q}(\sqrt{-3})$ can be a subfield of the invariant trace field.  This implies that the knots do not admit hidden symmetries via the Neumann and Reid characterization discussed above. 
 
In \cite{BBoWa}, conjecture \ref{3conj} is proven in the ``generic" case.   Namely, 

\begin{theorem} \cite{BBoWa} If $K$ is a hyperbolic knot whose complement does not admit hidden symmetries then $|C(K)| \leq 3$.
\end{theorem}
 
To prove this, it is shown that in the case when $K$ and $K'$ are commensurable and do not admit hidden symmetries, then the knot complements cyclically cover a common orbifold (as in Walter Neumann's examples). This generalizes the situation where a knot complement cyclically covers another knot complement, as in the case of Berge knots.  An {\it orbi-lens space} is the quotient of $S^3$ by a finite cyclic group.  The orbifold which is cyclically covered by the two knot complements is the complement of a knot in an orbi-lens space.   As another example of an orbifold, we give an example of an orbi-lens space.  Here the ramification locus is the cores of a genus 1 Heegaard splitting of the underlying space, which is a lens space. 
 
 \begin{description} 
 \item[5]  We consider $S^3$ as the unit 3-sphere in $\mathbb{C}^2$.  Then let $G$ be the group of isometries of $\mathbb{C}^2$ generated by $\phi$ where $\phi(z,w) = (e^{\frac{2 \pi i }{6}}z, e^{\frac{4 \pi i}{15}}w)$.  Then $G$ is cyclic of order 30 and leaves the unit three-sphere invariant.  $\phi^6$ fixes the $z$ axis of $S^3$, $\phi^{15}$ fixes the $w$ axis and $\phi^{10}$ acts freely.  Thus the quotient orbifold $\mathcal{L}$ has underlying space $|\mathcal{L}|$ a lens space with fundamental group of order $3$.  The ramification locus is two circles labeled 5 and 2 which are the cores of a genus 1 Heegaard splitting of  $|\mathcal{L}|$. 
 \end{description}

We conclude by remarking that understanding commensurability and hidden symmetries can lead to information about the symmetries of hyperbolic knot complements.  For example, we have the following:

\begin{lemma} \cite{BBoWa} If a hyperbolic knot $K$ does not admit hidden symmetries and $|C(K)| >1$, then $K$ is not amphichiral. If in addition $K$ is periodic, then it must be strongly invertible. 
\end{lemma}

{\bf Acknowledgements}  I learned much of this material from working with my coauthors Alan Reid, Steven Boyer and Michel Boileau and I thank them for their efforts.  I am also very grateful for the hospitality of the Harvard mathematics department and for support from the NSF.


\begin{thebibliography}{9999}
  
\bibitem{AR} I. R. Aitchison and J. H. Rubinstein, {\em Combinatorial
cubings, cusps, and the dodecahedral knots}, in
Topology '90, Ohio State Univ. Math. Res. Inst. Publ. {\bf 1}, 
pp. 17-26,  de Gruyter (1992).\\

\bibitem{BBoWa} M. Boileau, S. Boyer and G. S. Walsh, {\em On commensurability of knot complements}, Preprint. Available at http://www.tufts.edu/~gwalsh01/. \\


\bibitem{orbifoldbook} D. Cooper and C.  Hodgson and S. Kerckhoff, {\em Three-dimensional Orbifolds and Cone-Manifolds}, Math. Soc. of Japan, MSJ Memoirs {\bf 5} (2000). \\

\bibitem{CGLS} M. Culler, C. McA. Gordon, J. Luecke and P.B. Shalen, {\em Dehn surgery on knots}, Ann. of Math. {\bf 125} (1987), pp. 237--300.\\

\bibitem{GAW} F. Gonz{\'a}lez-Acu{\~n}a and W. C. Whitten, {\em Imbeddings of three-manifold groups}, Mem. Amer. Math. Soc.  {\bf 474} (1992).\\


\bibitem{GHH} 
O. Goodman, D. Heard and C. Hodgson,
{\em Commensurators of cusped hyperbolic manifolds}, Exp. Math. {17} (2008), pp. 283 -- 306. \\

\bibitem{FW}  B. Farb and S. Weinberger, {\em Hidden symmetries and arithmetic manifolds} in {\em  Geometry, spectral theory, groups, and dynamics}, Contemp. Math., {\bf 387} (2005), pp. 111--119. \\


\bibitem{Hof} N. Hoffman, preprint. \\

\bibitem{LLR} C. Leininger, D. D. Long and A. W. Reid, {\em Commensurators of non-free finitely generated Kleinian groups}, preprint. 

\bibitem{MM} M. L.  Macasieb and T. W. Mattman, {\em Commensurability classes of (-2,3,n) pretzel knot complements}, Algebraic and Geometric Topology {\bf 8} (2008) 1833--1853. \\

\bibitem{MR} C. Maclachlan and A. W. Reid, {\em The Arithmetic of Hyperbolic 3-manifolds}, Graduate Texts in Mathematics {\bf 219}, Springer-Verlag (2003).\\

\bibitem{outercircles} A. Marden, {\em Outer Circles: An Introduction to Hyperbolic 3-manifolds}, Cambridge University Press (2007). \\

\bibitem{Ma} G. Margulis, {\em Discrete Subgroups of Semi-simple Lie
Groups}, Ergeb. der Math. {\bf 17} Springer-Verlag (1989).\\

\bibitem{NR} W. D. Neumann and A. W. Reid, {\em Arithmetic of hyperbolic manifolds}, in
Topology '90, Ohio State Univ. Math. Res. Inst. Publ. {\bf 1} 
273--310,  de Gruyter (1992).\\


\bibitem{Re} A. W. Reid, {\em A note on trace-fields of Kleinian groups}, Bull. London Math. Soc.
{\bf 22} (1990), 349--352.\\

\bibitem{Re1}  A. W. Reid, {\em Arithmeticity of knot complements}, J. London Math. Soc. {\bf 43}
(1991), 171--184.\\

\bibitem{RW} A. W. Reid and G. S. Walsh, {\em Commensurability classes of two-bridge knot complements},   Algebraic and Geometric Topology, {\bf 8} (2008) 1031 -- 1057. \\


\bibitem{Ri1}  R. Riley, {\em Parabolic representations of knot groups. {I}},
Proc. London Math. Soc. {\bf 24} (1972), 217--242.\\

\bibitem{Ri2}  R. Riley, {\em Seven excellent knots}, Low-dimensional topology (Bangor, 1979),
L.M.S Lecture Note Series {\bf 48}, 81--151 C.U.P.(1982).\\

\bibitem{Sch} R. E. Schwartz, {\em The quasi-isometry classification of rank one lattices}, 
Publ. I.H.E.S {\bf 82} (1995), 133--168.\\

\bibitem{thurstonnotes} W. Thurston, {\em The Geometry and Topology of 3-manifolds} Princeton University lecture notes, 1980. Electronic version 1.1: http://www.msri.org/publications/books/gt3m/.  \\

\bibitem{WW} S. Wang and Y-Q. Wu, {\em Any knot complement covers
at most one knot complement}, Pacific J. Math. {\bf 158} (1993), 
387--395.\\

\end{thebibliography}
\end{document}